\documentclass[reqno,11pt]{amsart}
\usepackage{amssymb,amsmath}

\newtheorem{theorem}{Theorem}
\newtheorem{proposition}{Proposition}
\newtheorem*{conjecture*}{Conjecture}
\newtheorem*{corollary*}{Corollary}
\newtheorem{lemma}{Lemma}
\newtheorem{corollary}{Corollary}

\theoremstyle{definition}
\newtheorem{example}{Example}
\newtheorem{remark}{Remark}

\newcommand{\F}{{\mathbb F}}
\newcommand{\roots}{\varrho}
\newcommand{\st}{\;|\;}

\title{Complete Padovan sequences in finite fields}

\begin{document}

\maketitle

\begin{center}
{JUAN B. GIL} \\
{\small Penn State Altoona, 3000 Ivyside Park, Altoona, PA 16601}\\[1em]
{MICHAEL D. WEINER} \\
{\small Penn State Altoona, 3000 Ivyside Park, Altoona, PA 16601}\\[1em]
{CATALIN ZARA} \\
{\small Penn State Altoona, 3000 Ivyside Park, Altoona, PA 16601}
\end{center}

\bigskip

\section*{Abstract}
Given a prime $p\ge 5$, and given $1<\kappa<p-1$, 
we call a sequence $(a_n)_{n}$ in $\mathbb{F}_p$ 
a $\Phi_{\kappa}$-sequence if it is periodic with period $p-1$, and if it
satisfies the linear recurrence $a_n+a_{n+1}=a_{n+\kappa}$ with $a_0=1$. 
Such a sequence is said to be a complete $\Phi_{\kappa}$-sequence if in 
addition $\{a_0,a_1,\dots,a_{p-2}\}=\{1,\dots,p-1\}$. For instance, every 
primitive root $b$ mod $p$ generates a complete $\Phi_{\kappa}$-sequence 
$a_n=b^n$ for some (unique) $\kappa$.  A natural question is whether 
every complete $\Phi_{\kappa}$-sequence is necessarily defined by a 
primitive root.  For $\kappa=2$ the answer is known to be positive. 
In this paper we reexamine that case and investigate the case $\kappa=3$ 
together with the associated cases $\kappa=p-2$ and $\kappa=p-3$.

\section{Introduction}

For a prime number $p \ge 5$ and a number $\kappa\in \{2,\ldots,p-2\}$, 
a sequence $(a_n)_{n \in \mathbb{Z}}$ of elements of $\F_p$ is said to 
be a \emph{$\Phi_{\kappa}$-sequence} if
\begin{equation}\label{eq:Phi-seq}
 a_0=1 \quad \text{and} \quad a_{n+\kappa} = a_{n} + a_{n+1}
 \;\text{ for all } n \in \mathbb{Z},
\end{equation}
where ``='' means (throughout this paper) equality in $\F_p$. 
A $\Phi_{\kappa}$-sequence is called \emph{complete} if
\begin{gather}
\label{eq:Periodic}
 (a_n)_n \;\text{ is periodic, with period } p-1, \text{ and} \\ 
\label{eq:Permutation}
 \{ a_1, \ldots , a_{p-2} \} = \{ 2,\ldots, p-1 \}.
\end{gather}

The case $\kappa=2$ has been studied by Brison \cite{Bri92}. A 
$\Phi_{2}$-sequence satisfies a Fibonacci recurrence, so it is completely 
determined by the value of $a_1$. A complete $\Phi_{2}$-sequence is called 
a \emph{complete Fibonacci sequence}; for example, the only complete 
Fibonacci sequence in $\F_5$ is
$$\ldots, 4 ,2, 1, 3, 4,2,1,3,4, \ldots $$
with $a_0=1$ and $a_1=3$. Moreover, $x=3$ is a primitive root mod $5$, and 
it is a solution (in $\F_5$) of the equation $x^2=x+1$ (hence $x=3$ is a 
\emph{Fibonacci primitive root}), so the sequence can also be described as
$$ \ldots,3^{-2} (= 4), 3^{-1} (=2), 3^0 (=1), 
   3 ,3^2(=4), 3^3(=2), 3^4(=1), \ldots $$
The main result of \cite{Bri92}, stated here in a slightly weaker
form, generalizes this observation.
In Section~\ref{sec:Fibonacci} we will give an elementary proof of 
this theorem in order to motivate our approach.

\begin{theorem}[Brison] \label{th:Fibonacci} 
Let $p \geq 5$ be a prime number.  A $\Phi_{2}$-sequence $(a_n)_n$ is 
a complete Fibonacci sequence if and only if $a_n = b^n$ for all $n$, 
where $b$ is a Fibonacci primitive root.
\end{theorem}

The new results of this paper concern the case $\kappa=3$. Because
of the specific recurrence satisfied by $\Phi_3$-sequences
($a_{n+3} = a_{n+1} + a_n$), complete $\Phi_3$-sequences will be 
called \emph{complete Padovan sequences} \cite{Padovan}. Similar to the 
case $\kappa=2$, we will say that a primitive root $b$ in $\F_p$ is a 
\emph{Padovan primitive root} if it satisfies the condition $b^3=b+1$. 
Note that if $b$ is a Padovan primitive root, then the sequence
$$ \ldots ,b^{-1}(=b^{p-2}), 1, b, b^2, b^3, b^4, 
   \ldots, b^{p-2}, b^{p-1} (= 1), \ldots $$
is a complete Padovan sequence.  A natural question is whether these are 
the only examples of complete Padovan sequences in $\F_p$. Our results
state that this is the case, at least for certain prime numbers.

Let $\roots_p$ be the number of distinct roots of $X^3-X-1$ in $\F_p$.
The main results of this paper are the following two theorems.

\begin{theorem}\label{th:PadovanA} 
Let $p \geq 5$ be a prime number such that $\roots_p<3$.  
A $\Phi_{3}$-sequence $(a_n)_n$ is a complete Padovan sequence if and 
only if $a_n=b^n$ for all $n$, where $b$ is a Padovan primitive root.
\end{theorem}

In the case when $\roots_p=3$, we denote by $\alpha$, $\beta$, and 
$\gamma$ the roots of $X^3-X-1$ in $\F_p$. Further we let 
\begin{equation}\label{eq:orderN}
  N_p=\min\{|\alpha/\beta|,|\beta/\gamma|,|\gamma/\alpha|\}.
\end{equation} 

\begin{theorem}\label{th:PadovanB} 
Let $p \geq 5$ be a prime number such that $\varrho_p=3$ and 
$p\le N_p^2+1$. A $\Phi_{3}$-sequence $(a_n)_n$ is a complete Padovan
sequence if and only if $a_n=b^n$ for all $n$, where $b$ is a
Padovan primitive root.
\end{theorem}

We strongly believe that this theorem holds even if $p>N_p^2+1$. 
In fact, in Section~\ref{sec:powers} we will see that our condition on
$p$ can be relaxed, see \eqref{eq:weakCondition}. Numerical computations 
show that among all primes less than $10^5$, there are only 4 numbers that
cannot be handled by our proof of Theorem~\ref{th:PadovanB}, cf. 
Section~\ref{sec:examples}.  Nonetheless, for these cases one can 
manually check that the statement of our theorem is still true. 

In contrast to the Fibonacci case, the Padovan recurrence is of
order three, so in addition to $a_0=1$, one needs values for both
$a_1$ and $a_2$ to completely determine a $\Phi_3$-sequence. It is
therefore rather surprising that complete Padovan sequences are
determined by only one parameter. If $(a_n)_n$ is a complete
Padovan sequence, one can use an approach similar to the Fibonacci
case to get one condition relating $a_1$ and $a_2$. The difficulty
resides in proving a \emph{second} relation.

The ultimate goal will be to connect $\Phi_{\kappa}$-sequences in 
$\F_p$ with primitive roots of $p$. It is easy to see that if $p \geq 5$ 
is a prime, and $b \in \F_p$ is a primitive root mod $p$, then there 
exists a unique value $\kappa \in \{2, 3, \ldots , p-2\}$ such that 
$b^{\kappa}=b+1$. Therefore the sequence
$$ \ldots, b^{-1} (=b^{p-2}), 1, b, b^2, \ldots, b^{\kappa-1},
b^{\kappa}, \ldots, b^{p-2}, b^{p-1} (=1), \ldots $$
is a complete $\Phi_{\kappa}$-sequence. Data collected so far
suggests that these are in fact the only complete
$\Phi_{\kappa}$-sequences. For a fixed value $\kappa \in \{2,
\ldots, p-2\}$, we say that a primitive root $b$ in $\F_p$ is a
\emph{$\Phi_{\kappa}$-primitive root} if $b^{\kappa} = b+1$.

\begin{conjecture*}
Let $p \geq 5$ be a prime number.  A $\Phi_{\kappa}$-sequence $(a_n)_n$ 
is complete if and only if $a_n=b^n$ for all $n$, where $b$ is a 
$\Phi_{\kappa}$-primitive root.
\end{conjecture*}

At the end of the paper we briefly discuss the relation between the
conjugate cases $\kappa$ and $p-\kappa$ and prove the statement for 
$\kappa=p-2$, $\kappa=p-3$, $\kappa=\frac{p-1}{2}$, and 
$\kappa=\frac{p+1}{2}$.

\section{Fibonacci primitive roots}
\label{sec:Fibonacci}

In this section we discuss the characterization of complete Fibonacci 
sequences in terms of Fibonacci primitive roots, and give an elementary 
proof of Theorem~\ref{th:Fibonacci}. The key argument is the same as in 
\cite{Bri92}, but our approach is more direct.

\begin{proof}[Proof of Theorem~\ref{th:Fibonacci}]
It is not hard to see that a Fibonacci primitive root generates a
complete Fibonacci sequence: If $b$ is a Fibonacci primitive root
and $a_1=b$, then $a_2=a_0+a_1 = 1+b = b^2$, and, by induction,
$a_m = b^m$ for all integers $m$. Since $b$ is a \emph{primitive}
root, it follows that the sequence $(a_n)_n$ satisfies both the
periodicity and the completeness conditions, hence it is a
complete Fibonacci sequence.

The less trivial part is to show that for every complete Fibonacci
sequence, $a_1$ is a Fibonacci primitive root. The case $p=5$ can be 
checked separately, so from now on we assume that $p \geq 7$. Let 
$(a_n)_n$ be a complete Fibonacci sequence in $\F_p$ and let $b=a_1$. Let
\begin{equation}\label{eq:genseries}
 P(X) = \sum_{n=0}^{p-2} a_n X^n \in \F_p[X].
\end{equation}
Then the recurrence of $a_n$ implies
\begin{equation}\label{eq:polyident} 
  (1-X-X^2) P(X) = \bigl( 1-X^{p-1} \bigr) \bigl( 1+ a_{p-2}X \bigr).
\end{equation}
The right-hand side of \eqref{eq:polyident} is identically zero on
$\F_p^*$, while $P(X)$ can have at most $p-2$ roots in $\F_p^*$. 
Therefore, $1-X-X^2$ has at least one root in $\F_p^*$, and thus it 
has both roots in $\F_p^*$.

Let $\alpha$ and $\beta$ be the solutions of $x^2=x+1$ in $\F_p$.
Then $\alpha \neq \beta$ (since $p \neq 5$), hence
\begin{equation*}
 a_n = A \alpha^n + B \beta^n
\end{equation*}
for all integers $n$, where
\begin{equation*}
 A= \frac{b-\beta}{\alpha - \beta} \quad \text{and} \quad 
 B= \frac{b-\alpha}{\beta - \alpha}\,.
\end{equation*}

If $k$ is an integer such that $1 \leq k \leq p-2$, then
$$ \sum_{n=0}^{p-2} a_n^k = \sum_{n=1}^{p-1} n^k = 0 .$$
Therefore,
\begin{equation}
\label{eq:AB}
  0 = \sum_{n=0}^{p-2} a_n^k 
  = \sum_{n=0}^{p-2} \bigl( A \alpha^n + B \beta^n \bigr)^k 
  = \sum_{j=0}^k \tbinom{k}{j} A^jB^{k-j} 
  \sum_{n=0}^{p-2}(\alpha^j \beta^{k-j})^n \,.
\end{equation}
However, if $x \neq 0$, then
$$ \sum_{n=0}^{p-2} x^n = 
\begin{cases}
    0, & \text{if } x \neq 1 \\
    p-1, & \text{if } x=1 \\
\end{cases}\;.
$$

The key ingredient in the proof is finding a value of $k$ for
which the last sum in \eqref{eq:AB} is zero for all but one $j$. 
Using $\alpha+\beta=1$ and $\alpha\beta =-1$, we see that for all 
primes $p \geq 7$, the smallest such value is $k=4$. In fact, we get
\begin{equation}
\label{eq:AB2}
 0 = \sum_{n=0}^{p-2} a_n^4 = -6A^2B^2,
\end{equation}
which implies that one of $A$ or $B$ is zero. Without loss of generality, 
assume $A=0$. Then $b=\beta$, so $b^2=b+1$. An inductive argument shows 
that $a_n = b^n$ for all integers $n$, and since $(a_n)_n$ is a complete 
Fibonacci sequence, $b$ must be a primitive root, hence $a_1 = b$ is a 
Fibonacci primitive root.
\end{proof}

\section{Complete Padovan sequences and primitive roots}
\label{sec:Padovan}

Let $(a_n)_n$ be a complete Padovan sequence. Note that the
periodicity \eqref{eq:Periodic} and the recurrence relation 
\eqref{eq:Phi-seq} (with $\kappa=3$) imply, as in \eqref{eq:polyident}, 
\begin{equation*}
 (1-X^2-X^3)P(X) =0 \;\text{ on } \F_p^*.
\end{equation*}
Consequently, there must be at least one element $r\in\F_p$ that solves 
the equation $1-x^2-x^3=0$. Thus $1/r\in\F_p$ is a root of $f(X)=X^3-X-1$. 

\begin{proof}[Proof of Theorem~\ref{th:PadovanA}]
We only need to prove that a complete Padovan sequence with initial
values $(a_0,a_1,a_2)=(1,b,c)$ gives rise to a primitive root $b\in\F_p$
with $b^2=c$ and $b^3=b+1$.

\medskip
\textit{Case 1: $f$ has exactly one root in $\F_p$.}
\smallskip

In this case, $1-X^2-X^3$ must have exactly one root $r\in\F_p$. Thus
the other $p-2$ nonzero elements of $\F_p$ are roots of $P(X)$ from
\eqref{eq:genseries}, so
\begin{equation}\label{eq:factorization}
  P(X)= a_{p-2}\prod_{\substack{i=1\\i\not=r}}^{p-1}(X-i).
\end{equation}
Note that the periodicity of $(a_n)_n$ implies
$(a_{p-1},a_p,a_{p+1})=(1,b,c)$, and by \eqref{eq:Phi-seq} we have
$a_{p-2}=c-1$ and $a_{p-3}=b-c+1$. Comparing the constant term and the
coefficient of $X^{p-3}$ in \eqref{eq:factorization}, we conclude that
\begin{equation*}
  1= a_{p-2}/r \quad\text{and}\quad a_{p-3}= a_{p-2}r,
\end{equation*}
and consequently,
\begin{equation*}
 c=r+1 \quad\text{and}\quad b=r^2+r.
\end{equation*}
Using that $1-r^2-r^3=0$ it is easy to check that $b^2=c$ and
$b^3=b+1$. Thus $a_n=b^n$, and because of \eqref{eq:Permutation},
$b$ is a primitive root.

\medskip
\textit{Case 2: $f$ has only two distinct roots in $\F_p$.}
\smallskip

The discriminant of $f(X)=X^3-X-1$ is $-23$, so the only way for $f$ 
to have only two distinct roots in $\F_p$ is when $p=23$. We use
the periodicity of $(a_n)_n$ to conclude that
$1=a_0=a_{p-1}=a_{22}=22c+13b+17$ and so $c=13b+16$ in $\F_{23}$.
Note that because of \eqref{eq:Permutation} we have
\begin{equation}\label{eq:CubicSum}
 \sum_{n=0}^{p-2} a_n^3 = \sum_{n=1}^{p-1} n^3 = 0
 \;\text{ in } \F_p.
\end{equation}
Using $c=13b+16$ we get the cubic equation
\begin{equation*}
 0 = \sum_{n=0}^{21} a_n^3 = 10b^3 + b^2 +20
\end{equation*}
whose solutions in $\F_{23}$ are $b=3$ and $b=10$. However, the sequence
generated by the initial values $(1,3,9)$ fails to be complete, so
$b=3$ is not an admissible choice for $b=a_1$. On the other hand, $b=10$ 
is indeed a primitive root mod $23$.
\end{proof}

\begin{proof}[Proof of Theorem~\ref{th:PadovanB}]
Let $\alpha, \beta, \gamma\in\F_p$ be the distinct roots of $f(X)$.
Note that in this case $p$ must be different from $23$, so
$2\alpha+3$, $2\beta+3$, and $2\gamma+3$ are all different from $0$.
Note also that these roots satisfy the equations
\begin{equation}\label{eq:roots}
 \alpha+\beta+\gamma=0,\quad \alpha\beta\gamma=1, \;\text{ and }\;
 \alpha\beta+\alpha\gamma+\beta\gamma=-1.
\end{equation}
The recurrence relation \eqref{eq:Phi-seq} with initial values
$(1,b,c)$ gives the formula
\begin{equation*}
 a_n=A\alpha^n + B\beta^n + C\gamma^n \;\text{ for every $n$},
\end{equation*}
where
\begin{equation}\label{eq:ABC-coeff}
 A=\frac{\alpha^2 b + \alpha c + 1}{2\alpha+3},\quad
 B=\frac{\beta^2 b + \beta c + 1}{2\beta+3},\quad
 C=\frac{\gamma^2 b + \gamma c + 1}{2\gamma+3}.
\end{equation}

We make again use of \eqref{eq:CubicSum} to get
\begin{align*}
 0 &= \sum_{n=0}^{p-2} (A\alpha^n + B\beta^n + C\gamma^n)^3 \\
 &= \sum_{n=0}^{p-2}\; \sum_{i+j+k=3} \frac{3!}{i!j!k!}
 A^iB^jC^k(\alpha^{i}\beta^{j}\gamma^{k})^n \\
 &= \sum_{i+j+k=3} \frac{3!}{i!j!k!} A^iB^jC^k
    \sum_{n=0}^{p-2}(\alpha^{i}\beta^{j}\gamma^{k})^n.
\end{align*}
Now, using \eqref{eq:roots} it can be easily seen that
$\alpha^{i}\beta^{j}\gamma^{k}\not=1$ unless $i=j=k=1$.
Thus, for $(i,j,k)\not=(1,1,1)$, we get
\begin{equation*}
 \sum_{n=0}^{p-2}(\alpha^{i}\beta^{j}\gamma^{k})^n=
 \frac{1-(\alpha^{i}\beta^{j}\gamma^{k})^{p-1}}%
 {1-\alpha^{i}\beta^{j}\gamma^{k}}=0\; \text{ in } \F_p.
\end{equation*}
Finally, we arrive at the identity
\begin{equation}\label{eq:ABC}
 0 = \sum_{n=0}^{p-2} a_n^3 = -6 ABC
\end{equation}
which implies that at least one of the factors must vanish, say $C=0$,
and so the closed form of $a_n$ reduces to
\begin{equation}\label{eq:ClosedForm}
 a_n=A\alpha^n + B\beta^n \;\text{ for every $n$}.
\end{equation}

Note that $C=0$ implies $c=-\gamma b -\frac1{\gamma}$. If we can prove
that $b=\alpha$ or $b=\beta$, then this condition on $c$ together with
\eqref{eq:roots} give the desired identities $b^2=c$ and $b^3=b+1$ and
the theorem is proved.

Given the numbers $p$, $\alpha$, $\beta$ as above, for 
$k\in\{1,\dots,p-2\}$ we consider the set
$I_k=\{j\in\mathbb{Z}\st 0\le j\le k \text{ and }
\alpha^j\beta^{k-j}\equiv 1 \text{ mod } p\}$.
In the next section we will discuss some properties of this set and will
prove (Corollary~\ref{cor:ConditionN}) that for $p$, $\alpha$, $\beta$ 
as in this theorem, there is always a $k$, $1<k<p-1$, such that $I_k$ 
contains exactly one element.  If we let $k$ in \eqref{eq:AB} be such 
that $I_{k}=\{j_0\}$, then we get
\begin{equation}\label{eq:ABj0}
 0=\sum_{n=0}^{p-2} a_n^k = -\binom{k}{j_0} A^{j_0}B^{k-j_0}
\end{equation}
since the sum $\sum_{n=0}^{p-2}(\alpha^{j}\beta^{k-j})^n$ vanishes
in $\F_p$ for every $j$ with $\alpha^{j}\beta^{k-j}\not=1$. Therefore, 
either $A$ or $B$ must be zero. Now, together with the fact that $C=0$, 
the equations \eqref{eq:ABC-coeff} give $b=\alpha$ or $b=\beta$. 
This proves that $a_n=b^n$ for every $n$, and since $(a_n)_n$ is complete, 
$b$ is a Padovan primitive root.
\end{proof}

\section{Sum of powers and minimal exponent}
\label{sec:powers}

A crucial idea in the proof of Theorem~\ref{th:Fibonacci} and
Theorem~\ref{th:PadovanB} is to consider the sum of powers
$\sum_{n=0}^{p-2} a_n^k$ in $\F_p$ for some $k\in\{2,\dots,p-2\}$.
Since this sum is always zero, the aim is to find a suitable
exponent $k$ that allows us to extract useful information about the
sequence. For instance, for the Fibonacci sequence the exponent $k=4$
was a good choice. For a complete Padovan sequence, however, the situation 
is more subtle. In a first step, the choice $k=3$ allowed us to reduce
the closed form of $a_n$ to a sum of two terms, cf. \eqref{eq:ClosedForm},
but when working with the reduced form, the choice of $k$ is not clear at
all and depends on the prime number $p$ at hand.

Let $p$ be such that $X^3-X-1$ has three distinct roots
$\alpha, \beta, \gamma$ in $\F_p$. Let $(a_n)_n\subset \F_p$ be a
complete Padovan sequence with initial values $(1,b,c)$.
Assume  $\gamma^2 b +\gamma c+1=0$, i.e., $C=0$ in \eqref{eq:ABC-coeff}
so that $a_n$ reduces to \eqref{eq:ClosedForm}.
Under these assumptions we consider the set
\begin{align*}
I_k &=\{j\in\mathbb{Z}\st 0\le j\le k \text{ and }
\alpha^j\beta^{k-j}\equiv 1 \text{ mod } p\}
\intertext{and its dual}
I_k' &=\{j\in\mathbb{Z}\st 0\le j\le k \text{ and }
\alpha^{k-j}\beta^{j}\equiv 1 \text{ mod } p\}.
\end{align*}
Observe that $j\in I_k$ if and only if $k-j\in I_k'$ so that these sets
essentially contain the same information. 

Let $N=|\alpha/\beta|$ be the order of $\frac{\alpha}{\beta}$ in $\F_p$.
That is, $\big(\frac{\alpha}{\beta}\big)^N=1$ and
$\big(\frac{\alpha}{\beta}\big)^j\not=1$ for every
$j\in\{1,\dots,N-1\}$. In our situation it is easy to check that $N>3$.
Observe that
\begin{equation*}
1,\frac{\alpha}{\beta},\Big(\frac{\alpha}{\beta}\Big)^{2},\dots,
\Big(\frac{\alpha}{\beta}\Big)^{N-1}
\end{equation*}
are the (distinct) $N$th roots of unity.

\begin{lemma} \label{lem:kMinimal}
The order of $\alpha^N$ and $\beta^N$ in $\F_p$ is $(p-1)/N$.
\end{lemma}
\begin{proof}
Let $m\ge 1$ be such that $mN\le p-1$ and $\alpha^{mN}=1$.
Then $\beta^{mN}=1$ and so by \eqref{eq:ClosedForm} we must have
$a_{mN}=1$. This implies $mN=p-1$ so $m=(p-1)/N$.
\end{proof}

\begin{lemma}
$I_k\not=\varnothing$ for some $k\in\{1,\dots,p-2\}$.  Moreover, in this
case, we have $k=\ell(p-1)/N$ for some $\ell\in\{1,\dots,N-1\}$.
\end{lemma}
\begin{proof}
Suppose $I_k=\varnothing$ for every $k$. Then, in particular,
$\beta^k\not=1$ for $k=1,\dots,p-2$, so that $\beta$ must be a
primitive root mod $p$. Let $1<t<p-1$ be such that $\alpha=\beta^t$.
Thus $\alpha^j\beta^{k-j}=\beta^{(t-1)j+k}$ for every $j$ and $k$. If
we pick $k=p-t$ and $j=1$, then $\alpha\beta^{p-t-1}=\beta^{p-1}=1$
which implies $I_{p-t}\not=\varnothing$ and we get a contradiction.

Now let $k$ be such that $\alpha^j\beta^{k-j}=1$ for some $0\le j\le k$.
Then
\[ 1=(\alpha^{j}\beta^{(k-j)})^N
   =\Big(\frac{\alpha}{\beta}\Big)^{jN}\beta^{kN}=(\beta^{N})^k \]
which by Lemma~\ref{lem:kMinimal} implies that $(p-1)/N$ must divide $k$.
\end{proof}

\begin{lemma}\label{lem:kNcondition}
Let $I_k\not=\varnothing$ and let $j_0=\min(I_k)$. If $k< N+j_0$,
then $I_k=\{j_0\}$.
\end{lemma}
\begin{proof}
Let $j_1>j_0$ be such that $\alpha^{j_1}\beta^{k-j_1}=1$.
Thus $\alpha^{j_1}\beta^{k-j_1}=\alpha^{j_0}\beta^{k-j_0}$ and so
$\Big(\dfrac{\alpha}{\beta}\Big)^{j_1-j_0}=1$. But this implies
$j_1-j_0 =\ell N$ for some $\ell\ge 1$. Hence $j_1\ge N+j_0>k$ and
therefore $j_1\not\in I_k$.
\end{proof}

Let $k_{\min}$ denote the smallest $k$ for which $I_k\not=\varnothing$.

\begin{lemma}\label{lem:kmin}
If $k_{\min}>\frac{p-1}{N}$, then $k_{\min}< N+j_0$ and therefore 
$I_{k_{\min}}=\{j_0\}$.
\end{lemma}
\begin{proof}
Let $k_{\min}=\ell(p-1)/N$ for some $1<\ell<N$, so
\begin{equation}\label{NlPower}
 1=\alpha^{j_0}\beta^{k_{\min}-j_0}
 =\Big(\frac{\alpha}{\beta}\Big)^{j_0}\beta^{\ell(p-1)/N}
 =\Big(\frac{\alpha}{\beta}\Big)^{N+j_0}\beta^{\ell(p-1)/N}.
\end{equation}
Moreover, since $\beta^{(p-1)/N}$ is a $N$th root of unity,
we have $\beta^{(p-1)/N}=\Big(\dfrac{\alpha}{\beta}\Big)^{m}$ for some
$0\le m<N$. Therefore,
\begin{equation}\label{m(l-1)Power}
 1=\Big(\frac{\alpha}{\beta}\Big)^{j_0}\beta^{\ell(p-1)/N}
  =\Big(\frac{\alpha}{\beta}\Big)^{m+j_0}\beta^{(\ell-1)(p-1)/N}
\end{equation}
which in particular implies $m+j_0>(\ell-1)(p-1)/N$ since
$I_{(\ell-1)(p-1)/N}=\varnothing$.

Dividing \eqref{NlPower} by \eqref{m(l-1)Power} we get the equation
\begin{equation}\label{Contra}
 1=\Big(\frac{\alpha}{\beta}\Big)^{N-m}\beta^{(p-1)/N}.
\end{equation}
Now, if $k_{\min}=\ell(p-1)/N\ge N+j_0$, then we have
\begin{equation*}
  (\ell-1)(\tfrac{p-1}{N})<m+j_0<N+j_0\le \ell(\tfrac{p-1}{N})
\end{equation*}
which implies $N-m<(p-1)/N$. But the equation \eqref{Contra} would then
imply that $I_{(p-1)/N}\not=\varnothing$ contradicting the minimality of
$k_{\min}=\ell(p-1)/N$.
\end{proof}

\begin{corollary}\label{cor:ConditionN}
Let $N_p$ be as in \eqref{eq:orderN}. If $p\le N_p^2+1$, then 
$I_{k_{\min}}=\{j_0\}$.
\end{corollary}
\begin{proof}
By Lemma~\ref{lem:kmin} we only need to check the case when 
$k_{\min}=\frac{p-1}{N}$. Let $k=k_{\min}$. If $j_0=\min(I_k)$ and 
$j_0'=\min(I_k')$, then $(j_0,j_0')\not=(0,0)$. Otherwise it would 
imply $\alpha^k=\beta^k=1$ and so $a_k=A\alpha^k+B\beta^k=1$. 
But this contradicts the fact that, by definition, $a_k\not=1$ for 
$0<k<p-1$. Thus we can assume $j_0>0$. Then
\[ p\le N_p^2+1 \;\Rightarrow\; p-1\le N^2 \;\Rightarrow\; 
   k_{\min}=\frac{p-1}{N}\le N < N+j_0. \]
The statement now follows from Lemma~\ref{lem:kNcondition}.
\end{proof}

\begin{remark}\label{rem:conditionN}
According to our previous discussion, it is evident that the condition 
$p\le N_p^2+1$ in Theorem~\ref{th:PadovanB} can be replaced by the weaker 
condition $\frac{p-1}{N}<N+j_0$, or equivalently, 
\begin{equation}\label{eq:weakCondition}
 p< N^2+j_0 N+1. 
\end{equation}
Observe that if $k_{\min}>\frac{p-1}{N}$, then \eqref{eq:weakCondition} 
is automatically satisfied by Lemma~\ref{lem:kmin}. Thus we only need 
to  request \eqref{eq:weakCondition} for the cases when $k_{\min}=
\frac{p-1}{N}$. Some examples will be discussed in the next section.
\end{remark}

\section{Examples and further remarks}
\label{sec:examples}

Throughout this section we let $f(X)=X^3-X-1 \in \F_p[X]$.

\begin{example}
The set of numbers
\begin{multline*} 
\{7,11,17,37,67,83,113,199,227,241,251,271,283,367,373, 401,433,457,\\
  479,569,571,593,613,643,659,701,727,743,757,769,839,919,941,977\} 
\end{multline*}
contains all prime numbers $<1000$ for which there is a complete 
Padovan sequence  in $\F_p$ and $f(X)$ has exactly one root. 
For instance, for $p=7$ this root is $b=5$. It can be easily 
checked that $(a_0,a_1,a_2)=(1,5,4)$ generates a complete Padovan 
sequence, and that any other choice of initial values will not give  
such a sequence. Of course, $b=3$ is also a primitive root mod 7 and
$a_n=3^n$ is a complete sequence, but $f(3)\not=0$. However, $b=3$ 
solves the equation $x^4-x-1=0$ in $\F_7$, so $3^n$ is a complete 
$\Phi_4$-sequence in $\F_7$. 
\end{example}

\begin{example}
The set of numbers 
\begin{multline*} 
\{ 59,101,167,173,211,271,307,317,593,599,\\
   607,691,719,809,821,829,853,877,883,991,997 \}
\end{multline*}
contains all prime numbers $<1000$ for which there is a complete
Padovan sequence in $\F_p$ and $f(X)$ has three distinct roots.  
With each number we can associate some data according 
to our discussion in the previous section. In the following table we show 
a few examples. Under ``roots'' we list pairs of roots of $f(X)$, say 
$\alpha, \beta\in\F_p$, that generate complete Padovan sequences.  
Recall that $N=|\alpha/\beta|$.    
\medskip
\begin{center}
\begin{tabular}{|c|c|c|c|c|c|c|} \hline
$p$ & roots & $N$ & $k_{\min}$ & $\frac{p-1}{N}$ & $j_0$ & $j_0'$ \\\hline
59 & 13, 42 & 29 & 10 & 2 & 7 & 3 \\ \hline
101& 20, 89 & 20 & 20 & 5 & 16& 4 \\ \hline
101& 89, 93 & 25 & 8  & 4 & 7 & 1 \\ \hline
167& 134,73 & 83 & 14 & 2 & 5 & 9 \\ \hline
173& 97,110 & 86 & 10 & 2 & 1 & 9 \\ \hline
211& 205,97 & 15 & 14 & 14& 3 & 11\\ \hline
211& 97,120 & 42 & 30 & 5 & 6 & 24\\ \hline
271& 145,46 & 135& 22 & 2 & 17& 5 \\ \hline
307& 157,50 & 17 & 18 & 18& 11& 7 \\ \hline
307& 50,100 & 102& 15 & 3 & 4 & 11\\ \hline
307& 100,157& 102& 45 & 3 & 3 & 42\\ \hline
\end{tabular}
\end{center}

\medskip
Note that for $p=307$ with $N=17$ we have $p>N^2+1$. Nonetheless, as
discussed in Remark~\ref{rem:conditionN}, Theorem~\ref{th:PadovanB} still 
holds since $N+j_0>\frac{p-1}{N}$.
\end{example}

The relation between $j_0$, $j_0'$, and $k_{\min}$ is not a coincidence.
As a matter of fact, since 
$(\frac{\alpha}{\beta})^{k_{\min}}=(\frac{\alpha}{\beta})^{j_0+j_0'}$,
we have
\[ k_{\min} = j_0+j_0'+\ell N \]
for some integer $\ell\ge 0$. If $\ell=0$, then $N+j_0>\frac{p-1}{N}$ and
our proof of Theorem~\ref{th:PadovanB} works. 

\begin{example}
The set of numbers
\[ \{307,5851,24697,34961,87623,98801\} \]
contains all prime numbers $<10^5$ for which there is a complete
Padovan sequence in $\F_p$, $f(X)$ has three distinct roots, $p>N^2+1$, 
and $k_{\min}=\frac{p-1}{N}$ (cf.  Remark \ref{rem:conditionN}). 
More precisely, we have
\medskip
\begin{center}
\begin{tabular}{|c|c|c|c|c|c|} \hline
$p$ & $N$ & $k_{\min}$ & $j_0$ & $j_0'$ & $\ell$ \\\hline
307   &17  &18   &11  &7   &0 \\ \hline
5851  &39  &150  &4   &29  &3 \\ \hline
24697 &63  &392  &59  &18  &5 \\ \hline
34961 &92  &380  &89  &15  &3 \\ \hline
87623 &227 &386  &175 &211 &0 \\ \hline
98801 &52  &1900 &47  &33  &35 \\ \hline
\end{tabular}
\end{center}

\medskip
As mentioned above, $p=307$ and $p=87623$ are covered by our current 
proof of Theorem~\ref{th:PadovanB}. For the other four cases, the 
statement of the theorem can be checked by hand (computer).  
\end{example}

\begin{conjecture*}
The statement of Theorem~\ref{th:PadovanB} is true even if $p>N^2+1$.
\end{conjecture*}

\begin{remark}
In the case when $f(X)$ has three distinct roots in $\F_p$, we showed in
the proof of Theorem~\ref{th:PadovanB} that a complete Padovan sequence 
$(a_n)_n$ with $a_n=A\alpha^n+B\beta^n+C\gamma^n$ reduces to 
$a_n=A\alpha^n+B\beta^n$ with $A$ and $B$ as in \eqref{eq:ABC-coeff}. 
Moreover, because of \eqref{eq:roots}, $a_n$ can be written as
\[ a_n=-\gamma a_{n-1}-\frac{1}{\gamma} a_{n-2}\;\text{ for every } n, \]
thus it is in fact a generalized Fibonacci sequence with characteristic
polynomial $g(t)=t^2+\gamma t+\frac{1}{\gamma}=(t-\alpha)(t-\beta)$. 
According to \cite{BriNog03}, the set $\{1,a_1,\dots,a_{p-2}\}$ is then 
a standard $g$-subgroup. 
\end{remark}

\begin{remark}
The first part of the proof of Theorem~\ref{th:Fibonacci} works in a 
more general context: Let $(a_n)_n$ be a periodic sequence in $\F_p$, 
with period $p-1$, and satisfying a linear recurrence of order $\kappa$
\begin{equation*}
a_{n+\kappa} = s_0 a_n + s_1a_{n+1} + \cdots + s_{\kappa -1} 
a_{n+\kappa -1} \;\text{ for all } n,
\end{equation*}
with $a_0=1$. If
$$ S(X) = 1-s_{\kappa-1}X - \cdots - s_1 X^{\kappa-1} - s_0 X^{\kappa}
   \in \F_p[X],$$
and, as before,
$$ P(X) = \sum_{n=0}^{p-2} a_n X^n = 1 + a_1X + \cdots + a_{p-2}X^{p-2} 
   \in\F_p[X],$$
then $S(X)P(X) = (1-X^{p-1}) Q(X)$ for some $Q(X)\in \F_p[X]$, so $S(X)$ 
has at least one root $r\in\F_p$. Since
$$ X^{\kappa} - s_{\kappa-1}X^{\kappa -1} - \cdots - s_1 X - s_0
   =X^{\kappa}S(1/X), $$
it follows that the characteristic polynomial of the linear
recurrence has at least one root $1/r\in\F_p$. When $\kappa=2$, 
this implies that both roots are in $\F_p$, but for $\kappa \geq 3$
(as we have seen for $\kappa=3$), the situation is more complicated. 
Moreover, even if the roots are distinct and are all in $\F_p$, finding 
the right value(s) of $k$ in order to get relation(s) of the form 
\eqref{eq:AB2}, \eqref{eq:ABC}, or \eqref{eq:ABj0} is not obvious.
\end{remark}

\begin{proposition}\label{prop:conjugate}
Let $p \geq 5$ be a prime number. A $\Phi_{\kappa}$-sequence $(a_n)_n$
is complete if and only if $(a_{p-1-n})_n$ is a complete 
$\Phi_{p-\kappa}$-sequence.
\end{proposition}

\begin{proof}
Let $(a_n)_n$ be a complete $\Phi_{\kappa}$-sequence. Let
\[ \tilde a_n=a_{p-1-n} \;\text{ for every } n. \]
By definition, $(\tilde a_n)_n$ satisfies \eqref{eq:Periodic} and
\eqref{eq:Permutation}. So we only need to check \eqref{eq:Phi-seq}:
\begin{align*}
 \tilde a_n + \tilde a_{n+1}
 &= a_{p-1-n} + a_{p-1-n-1} \\
 &= a_{p-1-n-1+\kappa} = a_{p-1-n-1+p-p+\kappa}\\ 
 &= a_{p-1-(n+p-\kappa)} = \tilde a_{n+p-\kappa}. 
\end{align*}
Thus $(\tilde a_n)_n$ is a complete $\Phi_{p-\kappa}$-sequence. 
\end{proof}

This proposition together with Theorem~\ref{th:Fibonacci},
Theorem~\ref{th:PadovanA}, and Theorem~\ref{th:PadovanB} give us the
following corollaries.

\begin{corollary*}
Let $p \geq 5$ be a prime number.  A $\Phi_{p-2}$-sequence $(a_n)_n$ is 
a complete $\Phi_{p-2}$-sequence if and only if $a_n = b^n$ for all $n$, 
where $b$ is a $\Phi_{p-2}$-primitive root.
\end{corollary*}

\begin{corollary*}
Let $p \geq 5$ be a prime number such that $\roots_p<3$.  
A $\Phi_{p-3}$-sequence $(a_n)_n$ is a complete $\Phi_{p-3}$-sequence if 
and only if $a_n=b^n$ for all $n$, where $b$ is a $\Phi_{p-3}$-primitive root.
\end{corollary*}

\begin{corollary*}
Let $p \geq 5$ be a prime number such that $\varrho_p=3$ and 
$p\le N_p^2+1$. A $\Phi_{p-3}$-sequence $(a_n)_n$ is a complete 
$\Phi_{p-3}$-sequence if and only if $a_n=b^n$ for all $n$, where $b$ is 
a $\Phi_{p-3}$-primitive root.
\end{corollary*}

We finish this section discussing the case $\kappa=\frac{p-1}{2}$.
The corresponding statement for the case $\kappa=\frac{p+1}{2}$ follows 
by means of Proposition~\ref{prop:conjugate}.

\begin{theorem}
Let $p \geq 5$ be a prime number.  A $\Phi_{\frac{p-1}{2}}$-sequence
$(a_n)_n$ is complete if and only if $a_n =b^n$ for all $n$, 
where $b$ is a $\Phi_{\frac{p-1}{2}}$-primitive root. Moreover, in this
case $b=p-2$.
\end{theorem}

\begin{proof}
Let $(a_n)_n$ be a complete $\Phi_{\frac{p-1}{2}}$-sequence, so
\[ a_n+a_{n+1}=a_{n+\frac{p-1}{2}} \;\text{ for every } n. \]
Then
\begin{align*}
 a_{n-1}+2a_n+a_{n+1} &= (a_{n-1}+a_n) + (a_n+a_{n+1}) \\
 &= a_{n+\frac{p-1}{2}-1} + a_{n+\frac{p-1}{2}}
  = a_{n-1+p-1} = a_{n-1},
\end{align*}
so $2a_n+a_{n+1}=0$ and therefore $a_{n+1}=-2a_n=(p-2)a_n$. But this
implies 
\begin{equation*}
 a_n=(p-2)^n \;\text{ for every } n.
\end{equation*}
Thus $p-2$ is a $\Phi_{\frac{p-1}{2}}$-primitive root since $(a_n)_n$ 
is complete.  
\end{proof}

\begin{corollary*}
Let $p \geq 5$ be a prime number.  A $\Phi_{\frac{p+1}{2}}$-sequence
$(a_n)_n$ is complete if and only if $a_n =b^n$ for all $n$, 
where $b$ is a $\Phi_{\frac{p+1}{2}}$-primitive root. Moreover, in this
case $b=\frac{p-1}{2}$.
\end{corollary*}


\bigskip
2000 Mathematics Subject Classification: 11B37, 11B50
\end{document}